\documentclass[11pt,a4paper]{amsart}

\usepackage{amsmath,amsthm,amsfonts,amssymb}
\usepackage[a4paper,left=25mm,right=25mm,top=30mm,bottom=30mm,marginpar=25mm]{geometry}
\usepackage{upref}
\usepackage{pgf}

\usepackage{graphicx}
\usepackage[latin1]{inputenc}
\numberwithin{equation}{section}

\newtheorem{theorem}{Theorem}[section]

\def\pc{^\mathrm{pc}}
\def\qc{^\mathrm{qc}}
\def\rc{^\mathrm{rc}}

\def\rank{\mathop{\mathrm{rank}}}

\def\R{\mathbb{R}}

\newcommand{\bref}[1]{(\ref{#1})}

\def\dist{{\mathrm{dist}}}
\def\diag{{\mathrm{diag}}}

\def\epsilon{\varepsilon}
\def\I{\mathbb{I}}

\def\threed{_{3d}}
\def\Wt{\widetilde W}
\def\tr{\mathrm{tr}}
\def\Uset{O}
\def\Usetcl{\overline O}
\def\hdet{\theta}
\def\O{{\rm O}}
\def\SO{{\rm SO}}

\newcommand{\barg}[1]{\bigl(#1\bigr)}
\newcommand{\Barg}[1]{\Bigl(#1\Bigr)}
\newcommand{\bset}[1]{\bigl\{#1\bigr\}}
\newcommand{\Bset}[1]{\Bigl\{#1\Bigr\}}

\newcommand{\bbar}[1]{\bigl|#1\bigr|}

\title[Relaxation of a two-well energy in two dimensions]
{Relaxation of a model energy for the cubic to tetragonal 
phase transformation in two dimensions}

\author{Sergio Conti}
\address{Institut f\"ur Angewandte Mathematik, Universit\"at Bonn, 
Endenicher Allee 60, 53115 Bonn, Germany.}
\email{sergio.conti@uni-bonn.de}

\author{Georg Dolzmann}
\address{Fakult\"at f\"ur Mathematik, Universit\"at Regensburg, 93040 Regensburg, Germany.}
\email{georg.dolzmann@mathematik.uni-r.de}

\begin{document}

\begin{abstract}
We consider a two-dimensional problem in nonlinear elasticity which
corresponds to the  cubic-to-tetragonal phase transformation. Our model
is frame invariant and the energy density 
is given  by the squared distance from two potential wells. 
We obtain the quasiconvex envelope of the energy density and therefore the
relaxation of the variational problem. Our result  includes  
the constraint  of positive determinant.
\end{abstract}

\subjclass[2010]{73C50,49J40,52A30}
\keywords{cubic to tetragonal phase transformation,
quasiconvexity, relaxation}

\thanks{This work was partially supported by the Deutsche Forschungsgemeinschaft
through the Forschergruppe 797  {\em ``Analysis and computation of
  microstructure in finite plasticity''}, project
 DO 633/2-2 (second author)
and Sonderforschungsbereich 1060 {\em ``The mathematics of emergent effects''}, project A6 (first author).}

\date{\today}

\maketitle

\section{Introduction and statement of the result}

Fully nonlinear models for the description of solid to solid phase
transformations within the mathematical framework of elasticity theory 
have attracted a lot of attention in the past twenty five years,
starting with the seminal
papers~\cite{BallJames87,BallJames92,ChipotKinderlehrer88}.
One of the questions of interest is the characterization 
of macroscopic models which capture the essential features of the 
mechanical behavior of a given system without resolving all 
the structures which may develop on small scales. In this context,
the theory of relaxation and
generalizations of the convex hulls of functions and sets play an important
role. The most important notion of convexity
is the notion of quasiconvexity in the sense of Morrey~\cite{Morrey1966}
and the related definitions of quasiconvex envelopes of functions and
hulls of sets~\cite{Dacorogna1989,RoubicekBook1997,MuellerLectureNotes}.

In this note we focus on a model for a specific solid to solid
phase transformation of austenite-martensite type, namely the 
cubic to tetragonal phase transformation. Suppose that the 
energy of such a material in the martensitic or high temperature phase
is characterized by an energy density $W\threed$ which we may assume to be
nonnegative with $K\threed=W\threed^{-1}(0)\neq \emptyset$. It follows from the 
invariance under change in observer and the symmetry of the underlying
point group that in the three-dimensional setting
\begin{align*}
 K\threed = \bigcup_{i=1}^3 \SO(3)U_i\,,\quad U_i = \frac{1}{\lambda}\I_3 + 
\barg{ \lambda^2-\frac{1}{\lambda} }e_i\otimes e_i\,,
\end{align*}
where  $\lambda\in (0,\infty)$, and where $(e_i)_{i=1, \ldots, n}$ 
denotes  for $n\in \mathbb{N}$
the canonical basis in $\R^n$, $\I_n$ the $n\times n$
identity matrix, and $\SO(n)$ the group of proper rotations
of $\R^n$. We chose 
$\det U_i=1$ since in most materials
the phase transformation occurs  without significant
change of volume. A basic model energy in this situation is given by
\begin{align}\label{3denergy}
 W\threed(F) = \dist^2\barg{F, K\threed} = \min_{i=1,2,3}\,
\min_{Q\in \SO(3)}\bbar{F-QU_i}^2\,.
\end{align}
To the best of our knowledge, both the computation of the full relaxation
$W\threed\qc$ of the energy density and the computation of the quasiconvex hull
$K\threed\qc$ of the set of its minimizers $K\threed$
are open problems. Partial results were obtained
in~\cite{Bhattacharya1992,DolzmannKirchheim2003,ContiDolzmannKirchheim2007}.

Considerable
progress has been achieved, however, in the case of two potential wells. The quasiconvex 
hull of two wells in three dimensions, i.e., of the set 
$\widetilde K\threed=\SO(3)U_1\cup \SO(3)U_2$, was computed in~\cite{BallJames92}, some
generalizations were analyzed in~\cite{DolzmannEtAl2000,Matos1999}, and 
various proofs were given for the characterization of the quasiconvex hull
of two wells in two dimensions, even in the more general case of 
two wells with different determinant, i.e., for the set
$K=\SO(2)\I_2 \cup \SO(2)\diag(\lambda,\mu)$ with $\lambda$, $\mu>0$,
see, e.g.,~\cite{Sverak1993,MuellerLectureNotes,DolzmannLN}. 
The linear case has been analyzed, also for more general energies,  in
\cite{KohnStrang1986,LurieChervaev1988,Kohn1991,Pipkin1991}.

In this paper we provide an explicit
relaxation formula for the analogue of~\bref{3denergy} in the corresponding
two-dimensional model, namely,
\begin{equation}\label{Wformula}
W(F)=
  \mathrm{dist}^2(F, K) + \hdet(\det F)\,,
\end{equation}
where $\hdet:\R\to[0,\infty]$ is convex and lower semicontinuous and
\begin{align}\label{Kformula}
K = \SO(2)U_1 \cup \SO(2)U_2\,,\quad \text{ where }\quad 
   U_1=
  \begin{pmatrix}
    \lambda&0\\0&1/\lambda
  \end{pmatrix}\,,\quad
  U_2=
  \begin{pmatrix}
    1/\lambda&0\\0&\lambda
  \end{pmatrix}\,,
\end{align}
for some fixed $\lambda>1$.
The function $\hdet$ permits to incorporate easily the constraint that 
the determinant of the deformation gradient $F$ has to be positive in order to
rule 
out interpenetration of matter. Here and in the following, $\dist(\cdot)$
denotes the Euclidean distance. 

One important ingredient in the proof is to view the energy in the appropriate
variables. Motivated by the formula for $\widetilde
K\threed\qc$ 
in~\cite{BallJames92} we set
\begin{align*}
  v=\frac1{\sqrt2}
  \begin{pmatrix}
    1\\1
  \end{pmatrix}\,, \hskip1cm
  w=\frac1{\sqrt2}
  \begin{pmatrix}
    1\\-1
  \end{pmatrix}
\end{align*}
and define for $F\in \R^{2\times 2}$ the coordinates
\begin{align}\label{coordinates}
  x(F)=|Fv| \,,\quad y(F)=|Fw| \,,\quad d(F)=\det F\,.
\end{align}
The main theorem can now be stated as follows.

\begin{theorem}\label{maintheorem}
Let $\lambda>1$, $K\subset\R^{2\times 2}$ be given by
(\ref{Kformula}),
let $\Uset=\bset{(x,y,d)\in\R^3,\,x,\,y>0,\,xy>|d|}$, define the functions 
$A$, $g:\Usetcl\to \R$ by
\begin{align}\label{Agfunctions}
\begin{split}
  A(x,y,d)&\, =(x^2+y^2)\frac{|U_1|^2}2 +
    (\lambda^2-\frac1{\lambda^2}) \sqrt{x^2y^2-d^2} + 2 d \,,\\
  g(x,y,d)&\, =x^2+y^2+|U_1|^2 -2 \sqrt{A(x,y,d)} \,,
\end{split}
\end{align}
and, for a given convex and lower semicontinuous function
$\hdet:\R\to[0,\infty]$, the energy density 
$W\colon \R^{2\times 2} 
\to \R\cup\{\infty\}$ by (\ref{Wformula}).

Then $W\rc = W\qc = W\pc$ and the quasiconvex envelope $W\qc$ is given by
  \begin{equation}\label{Wqcformula}
    W^{qc}(F) = 
h(x(F), y(F), d(F)) + \hdet(\det F)\,,
  \end{equation}
where $h: \Usetcl \to \R$ is defined by
\begin{align}\label{hfunction}
  h(x,y,d)&\, = \min_{\xi\in [x,\infty),\,\eta\in [y,\infty)} g(\xi, \eta, d)\,.
\end{align}

\end{theorem}
Remarks.
\begin{enumerate}
\item The minimum in the definition of $h$ can be computed explicitly in terms of
the roots of a polynomial equation of fourth order, see
Section~\ref{secformulaphasebdry}.
\item\label{enumincompr} As a special case, if $\hdet(t)=0$ for $t=1$, and $\infty$ otherwise, we
  obtain a result for incompressible materials.
\item The set $\Uset$ is not convex.
\item\label{enumorient}
A typical example in nonlinear elasticity is the choice of $\hdet\in C^1((0,\infty))$
which satisfies
$\lim_{t\to 0}\hdet(t) = \lim_{t\to\infty}\hdet(t) = \infty$, for example,
$\hdet(t) = \log^2(t)$, in order to obtain a density which rules out
interpenetration of  
matter (and $\hdet(t)=\infty$ for $t\le 0$).
\item If $\theta$  has at most linear growth, then 
$W$ has quadratic growth and 
the lower semicontinuous envelope of  $\int W(Du)dx$ is given by $\int W\qc(Du)dx$, see \cite{Dacorogna1989,MuellerLectureNotes}. 
A corresponding result incorporating constraints on the determinant of the type mentioned in \ref{enumincompr} and \ref{enumorient}
will be discussed elsewhere \cite{ContiDolzmannRelax}.  
\item 
For $\lambda=1$ our formula reduces to the
well-known relaxation of the squared  distance function to $\SO(2)$,
as given for example in \cite{Silhavy1999}.
\end{enumerate}

We recall that a function 
$f:\R^{2\times 2}\to \R\cup\{\infty\}$ is said to be quasiconvex if
$f(F)\le \int_{(0,1)^2} f(F+D\varphi)dx$ for all $F\in \R^{2\times 2}$ and $\varphi\in W^{1,\infty}_0((0,1)^2;\R^2)$. The function 
$f$ is said to be rank-one convex if
for all $F$, $R\in \R^{2\times 2}$ with $\rank(R)=1$ the function
$t\mapsto f(F+t R)$ is convex. The function $f$ is said to be polyconvex if 
it can be written as $f(F)=g(F, \det F)$ with 
$g:\R^{5}\to \R\cup\{\infty\}$ convex and lower semicontinuous. 
The rank-one convex, quasiconvex, and polyconvex
envelopes $f\rc$ and $f\pc$ are the largest rank-one convex, quasiconvex  and polyconvex
functions less than or equal to $f$, respectively. 
If $f$ is polyconvex, then it is also rank-one convex, hence $f\pc\leq f\rc$.
We refer to~\cite{Dacorogna1989,MuellerLectureNotes} for more information on these notions 
of convexity and their relations. In the following we use for 
two vectors $a$, $b\in\R^2$ and two matrices $A$, $B\in \R^{2\times 2}$ the
notation $a\cdot b$ and $A:B$ for the inner product in $\R^2$
and $\R^{2\times 2}$, respectively. Hence $|A|^2=A:A=\tr(AA^T)$. Finally
$a\otimes b \in \R^{2\times 2}$ is given by $(a\otimes b)_{ij}=a_i b_j$.
We write $\R_+=(0,\infty)$.

\section{Proof}\label{proofs}
The general strategy of the proof is to verify first that the rank-one
convex and the polyconvex envelope coincide. Let $\Wt$ be the 
formula on the right-hand side in~\bref{Wqcformula}. In order to prove 
$W\pc=W\rc=\Wt$ 
 it suffices to show that $W\rc\leq  \Wt\leq W$ and that 
$\Wt$ is polyconvex. Then $W\pc \leq W\rc\leq \Wt\leq  W\pc$ and therefore 
$W\rc=W\pc=\Wt$. The quasiconvex envelope will be discussed at the end. We divide
the proof into several steps. 

\medskip

\textit{Step~1: $\Wt\leq W$.} This follows immediately once we have
shown that 
  \begin{equation}\label{eqwlambdag}
    W(F) = 
g\barg{x(F), y(F), d(F)} + \hdet(\det F)\,.
  \end{equation}
In order to establish~\bref{eqwlambdag}, we denote 
the signed singular values of a matrix $F\in \R^{2\times 2}$ by
$\lambda_i=\lambda_i(F)$ and use 
the convention that $\lambda_2\ge|\lambda_1|\geq 0$, $\lambda_1\lambda_2=\det
F$. Since $|U_1|=|U_2|$ we obtain 
\begin{align}\label{criticalcalculation}
\begin{split}
\dist^2(F, \SO(2) U_1\cup \SO(2) U_2)&= \min_{i=1,2}\,  \min_{Q\in\SO(2)} |F-QU_i|^2
\\ &=
|F|^2 + |U_1|^2 - 2 \max_{i=1,2}\,\max_{Q\in\SO(2)} FU_i : Q\\
&= |F|^2 + |U_1|^2 - 2  \max_{i=1,2}(\lambda_2+\lambda_1)(FU_i)\\
&= |F|^2 + |U_1|^2 - 2   \max_{i=1,2}\sqrt{|FU_i|^2 + 2 \det F}\,.
\end{split}
\end{align}
Moreover,
\begin{align*}
  |FU_1|^2=\lambda^2 |Fe_1|^2+\frac1{\lambda^2} |Fe_2|^2=
 (\lambda^2+\frac1{\lambda^2}) \frac{|Fv|^2+|Fw|^2}{2}
 +  (\lambda^2-\frac1{\lambda^2}) Fv\cdot Fw\,.
\end{align*}
Since replacing $U_1$ by $U_2$ corresponds to interchanging $\lambda$ 
and $1/\lambda$, the maximum over $i$ is given by
\begin{align*}
  \max_{i=1,2} |FU_i|^2=
 (\lambda^2+\frac1{\lambda^2}) \frac{x^2(F)+y^2(F)}{2}
 +  (\lambda^2-\frac1{\lambda^2}) |Fv\cdot Fw|\,,
\end{align*}
where  we used the coordinates defined in~\bref{coordinates}.
From $\sin^2+\cos^2=1$ we obtain
\begin{align*}
  |Fv\cdot Fw|= \sqrt{x^2(F)y^2(F)-{\det}^2(F)}
\end{align*}
and therefore
\begin{align*}
  W(F)=x^2(F)+y^2(F)+|U_1|^2 -2 \sqrt{A(x(F),y(F),d)}+ \hdet(\det F)\,,
\end{align*}
where $A$ is given by~\bref{Agfunctions}. This establishes
 (\ref{eqwlambdag}).

\medskip

\textit{Step~2: The upper bound $W\rc\leq \Wt$.} This bound follows
from standard arguments based on a minimization along rank-one lines
which has been used in most examples for explicit relaxation results in
nonlinear elasticity and plasticity, see, 
e.g.,~\cite{AlbinContiDolzmann,CarstensenPlechac1997,ContiDolzmannKreisbeck2011,
ContiEtAlM3AS2013,ContiTheil2005,DeSimoneDolzmannARMA2002,
Silhavy1999,Silhavy2007} and the references therein.

In the following we use for $(x,y,d)\in \Usetcl$
the notation $z(x,y,d)=\sqrt{x^2y^2-d^2}$.
Note that $A$ and $g$ are continuous and nonnegative functions on $\Usetcl$.
Fix  a matrix $F$
and define the rank-one line $t\mapsto F_t$ by
\begin{align*}
  F_t=F(\mathrm{Id} + t v\otimes w)\,.
\end{align*}
Since the vectors $v$ and $w$ are orthogonal we obtain 
\begin{align*}
  d(F_t)=d(F)\,,\quad x(F_t)=x(F) \quad\text{ and }\quad y(F_t)=|Fw+ tFv|\,.
\end{align*}
If $Fv\ne0$ then
 $y(F_t)$ tends to infinity for $t\to \pm\infty$ and 
by continuity for every $y'>y(F)$ there exist $t_-<0<t_+$ such that
$y(F_{t_\pm})=y'$. Then there is $\mu\in(0,1)$ such that
 $F=\mu F_{t_+} + (1-\mu) F_{t_-}$,  and from the 
convexity of $W\rc$ along the rank-one line $t\mapsto F_t$ 
and (\ref{eqwlambdag}) we infer
\begin{align}\label{eqwrcupperb}
  W\rc(F)\le g(x(F), y', d(F))+ \hdet(\det F)\quad \text{ for all } y'\ge
  y(F)\,. 
\end{align}
If $Fv=0$ we take instead $F_t=F + t w\otimes w$. Then  $d(F_t)=d(F)=0$, 
$x(F_t)=x(F)=0$, $y(F_t)=|Fw+tw|$ tends to infinity for $t\to \pm\infty$
and  (\ref{eqwrcupperb}) follows as above.

If one interchanges $v$ and $w$ one obtains
\begin{align*}
  W\rc(F)\le h(x(F), y(F), d(F))+ \hdet(\det F)
\end{align*}
where $h$ is defined in~\bref{hfunction}. Note that by
 continuity and growth the infimum in (\ref{hfunction}) is attained.

\textit{Step~3: The lower bound $\Wt\leq W\pc$.}
The general idea is to show that $h$ is the restriction to $\Uset$ of a convex
function $f$, 
nondecreasing in
the first two arguments, which will be constructed below. Since $h$ is continuous on 
$\Usetcl$ the result 
will follow.

We compute for $(x, y, d)\in \Uset$  the gradient
\begin{align}\label{ggradient}
  D g(x,y,d) = D(x^2+y^2) - \frac{D A}{\sqrt A}(x,y,d)
\end{align}
and the matrix of second derivatives,
\begin{align*}
D^2 g = 2(e_1\otimes e_1+e_2\otimes e_2)  -  \frac{D^2 A}{\sqrt A}
+ \frac{D A\otimes DA}{2 A^{3/2}}\,.
\end{align*}
From the explicit form of $A$ we obtain with $z=\sqrt{x^2y^2-d^2}$
\begin{align*}
  \partial_x A(x,y,d) 
&\,=
x|U_1|^2 + (\lambda^2-\frac1{\lambda^2}) \frac{xy^2}{\sqrt{x^2y^2-d^2} }
=x|U_1|^2 + (\lambda^2-\frac1{\lambda^2}) \frac{xy^2}{z }
\end{align*}
and since $\partial_x A$ is linear in $x$, except from the dependence of $z$ on $x$,
\begin{align*}
  \partial_x \partial_x A(x,y,d) = \frac{\partial_x A(x,y,d)}{x} -
    (\lambda^2-\frac1{\lambda^2}) \frac{x^2y^4}{z^3}\,.
\end{align*}
The analogous formulas for the other derivatives imply that
\begin{align*}
  D^2A(x,y,d)&=
  \begin{pmatrix}
    \frac{\partial_x A(x,y,d)}{x} & 0 & 0\\
    0 & \frac{\partial_y A(x,y,d)}{y} & 0\\
    0 & 0 & 0
  \end{pmatrix}
 \\
& -    (\lambda^2-\frac1{\lambda^2})
 \frac{1}{z^3 }
  \begin{pmatrix}
x^2y^4  &  2d^2xy-x^3y^3 & -dxy^2\\
2d^2xy-x^3y^3&x^4y^2 & -dx^2y\\
-dxy^2 & -dx^2y & x^2y^2
  \end{pmatrix} \,.
\end{align*}
Note that 
\begin{align*}
 \partial_x g(x,y,d) = 2x - \frac{\partial_x A}{\sqrt{A}}(x,y,d)\geq 0 
\quad \Leftrightarrow\quad 
\frac{\partial_x g(x,y,d)}{x} = 
2 - \frac{\partial_x A(x,y,d)}{x\sqrt{A(x,y,d)}}\geq 0\,.
\end{align*}
Therefore
\begin{align*}
D^2 g(x,y,d) &=  \begin{pmatrix}
  \frac{\partial_x g(x,y,d)}{x} & 0 & 0\\
    0 & \frac{\partial_y g(x,y,d)}{y} & 0\\
    0 & 0 & 0
  \end{pmatrix}
+ \frac{D A\otimes DA}{2 A^{3/2}}(x,y,d) \\
& +    (\lambda^2-\frac1{\lambda^2})
 \frac{1}{z^3 \sqrt{A(x,y,d)}}
  \begin{pmatrix}
x^2y^4  &  2d^2xy-x^3y^3 & -dxy^2\\
2d^2xy-x^3y^3&x^4y^2 & -dx^2y\\
-dxy^2 & -dx^2y & x^2y^2
  \end{pmatrix} \,.
\end{align*}
The last matrix in this expression is positive semidefinite.
To see this, note that the determinant of the matrix
is equal to zero and that
the determinant of the first $2\times 2$ block is 
equal to $4d^2x^2 y^2(x^2y^2-d^2)$ and thus
positive on $\Uset$. 
The assertion follows now by continuity of the determinant if one adds 
$\epsilon>0$ to the $(3,3)$ element, computes the determinant, and considers
the limit as $\epsilon$ tends to zero.
Therefore each of the three terms is positive semidefinite in the region
\begin{align}\label{beginningchanges}
 V =\bset{ (x,y,d)\in \Uset,\,\partial_x g\ge 0\text{ and }\partial_y g\ge 0 }\,.
\end{align}
We conclude that $D^2g\ge 0$ in $V$.
To identify $V$, we study for fixed $y$ and $d$ the sign of the function 
\begin{align*}
  x\mapsto \partial_x g(x,y,d)\,,
\end{align*}
defined on  $x\in (|d|/y,\infty)$.
As $x\to |d|/y$ we have $z\to0$, hence
$\partial_x A\to\infty$ and $\partial_x g \to -\infty$
(if $d=0$, then still $\partial_xg/x\to-\infty$). Similarly, for
$x\to\infty$  we have $\partial_x g(x,y,d)\sim 2x + O(1)>0$.  Therefore there
is at least one zero in $x\in (|d|/y,\infty)$. The condition $\partial_x g=0$ is equivalent to 
\begin{align*}
  \frac{\partial_x A(x,y,d)}{x}=2\sqrt{A(x,y,d)}\,,
\end{align*}
which, in view of the foregoing definitions, is equivalent to
\begin{equation}\label{eqdxada}
|U_1|^2 +
    (\lambda^2-\frac1{\lambda^2}) \frac{y^2}{\sqrt{x^2y^2-d^2} }
=2\sqrt{(x^2+y^2)\frac{|U_1|^2}2 +
    (\lambda^2-\frac1{\lambda^2}) \sqrt{x^2y^2-d^2} + 2 d}\,.
\end{equation}
The left-hand side is strictly monotone decreasing in $x$, the right-hand side
strictly monotone increasing. Therefore there is at most one solution and since
we have already shown that there exists at least one solution, we conclude that  
for every pair $y>0$, $d\in\R$ there is a unique value $x=\varphi(y,d)\in(|d|/y,\infty)$ such that  
\begin{align*}
  \partial_x g(\varphi(y,d),y,d)=0\,.
\end{align*}
The above discussion then shows that
\begin{equation*}
  \partial_x g(\xi, y,d )\ge 0 \text{ for all } \xi\ge \varphi(y,d)\,.
\end{equation*}
By the implicit function theorem, the function
$\varphi:(0,\infty)\times\R\to(0,\infty)$ is smooth. By construction,
 $\varphi(y,d)>|d|/y$ and in particular $\varphi(y,d)$
tends to infinity as $y$ tends to zero. 
We show in Section \ref{secmonotonicity} below that $y\mapsto\varphi(y,d)$ is monotonically decreasing. Therefore $\varphi(y,d)$ has a finite limit
as $y$ tends to infinity. In particular, for all $d$
fixed the graph of $\varphi(\cdot, d)$ intersects the line $x=y$ and hence
there exists a unique $y$ such that $\varphi(y,d)=y$, 
i.e., 
\begin{align}\label{fixedpoint}
\text{ for all }d\text{ there exists a unique }
p>0\text{ with } p=\varphi(p,d)\,.
\end{align}
Interchanging 
$x$ and $y$ similar results, with the same function $\varphi$, hold for $\partial_yg$ and in particular $\partial_y g(x,\varphi(x,d),d)=0$.
Therefore
\begin{align*}
 V = \bset{ (x,y,d):
xy-|d|>0, \, x \geq \varphi(y,d),\,y\geq \varphi(x,d) }\,.
\end{align*}
To see how the two latter conditions interact, consider for definiteness the 
region  $\varphi(y,d)\le x\le y$. Then 
 $\partial_x g\ge 0$ or equivalently  $\partial_x A/x\le 2\sqrt A$. 
Since $x\le y$, from the explicit expression for $\partial_x A$ we obtain
$\partial_y A/y\le \partial_xA/x\le 2\sqrt A$, and therefore $\partial_y g\ge 0$ (all quantities are evaluated at $(x,y,d)$). We conclude that
\begin{equation}\label{eqfourthregion}
  \text{ if } \varphi(y,d)\le x\le y \text{ then } \varphi(x,d)\le y\,.
\end{equation}

\begin{figure}[t]
   \includegraphics[width=0.95\linewidth]{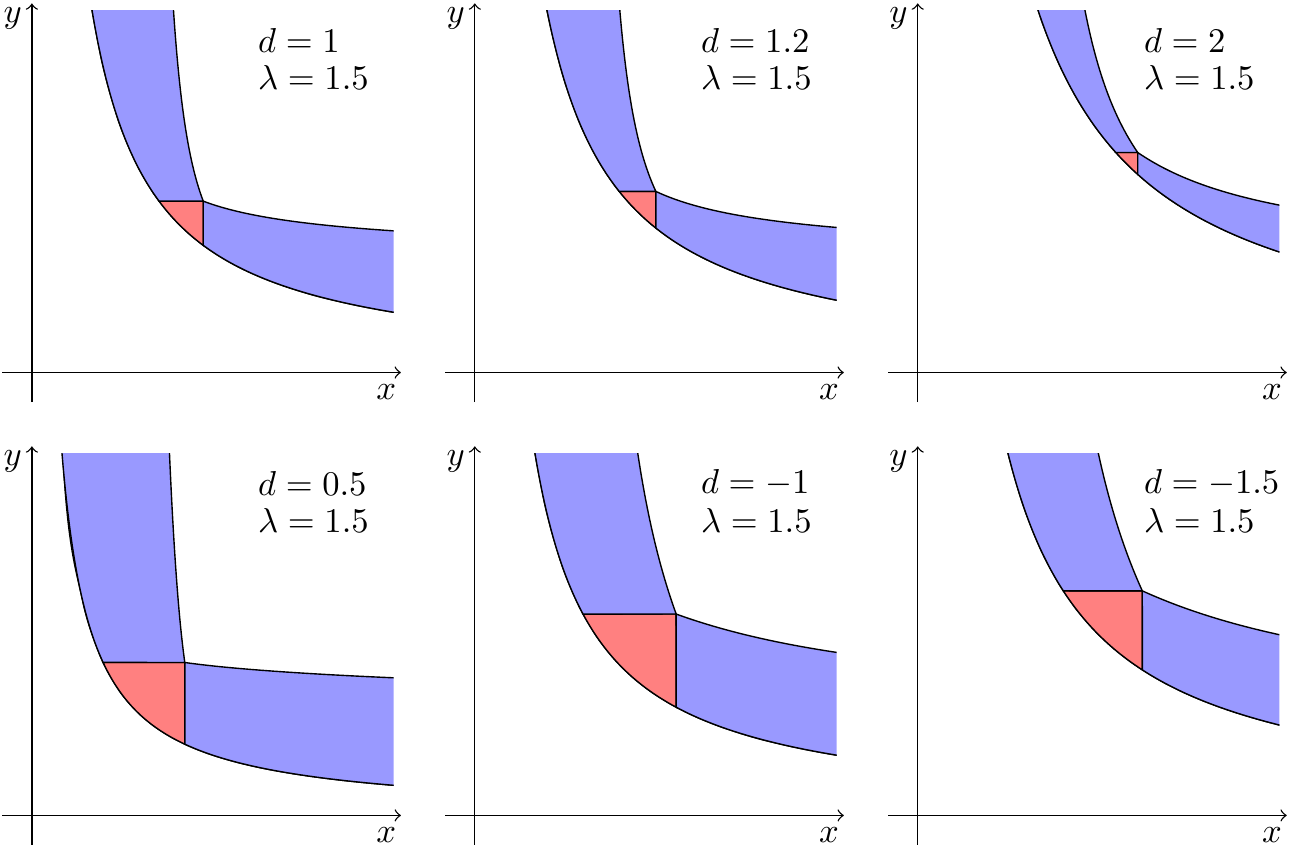}
\caption{Sketch of the phase diagram for $\lambda=1.5$, for various values of 
  the determinant $d$. The region below the
  hyperbola $xy=|d|$ is not accessible. The inner boundaries represent the
  curves  $x=\varphi(y,d)$ and $y=\varphi(x,d)$, which cross at $(p,p)$ as
  defined in (\ref{fixedpoint}). For $x<p$ and $y<p$ (red region) the energy
  is constant (for each $d$), in the 
  stripes close to the hyperbola ($|d|/y<x<\varphi(y,d)$ and the symmetric
  one, blue) the energy depends only on one of the two variables, in the 
  region above the curves (region $V$) it coincides with the unrelaxed energy. 
\label{fig1}}
\end{figure}

To show that $h(x(F), y(F), d(F))$ is polyconvex to show that 
$h$ is the restriction to $\overline\Uset$ of a function $f:\R^3\to \R$
which is convex on $\R^3$ and monotone increasing in its first two
arguments. 
Recall the definition of $p=p(d)$ in~\bref{fixedpoint}. We define $f$ by
\begin{align*}
  f(x,y,d)=
  \begin{cases}
    g(\varphi(y,d),y,d) & \text{ if } p\le y\text{ and } x\le \varphi(y,d)\,,\\
    g(p,p,d) & \text{ if } x\le p \text{ and } y\le p\,,\\
    g(x,\varphi(x,d),d) & \text{ if } p\le x\text{ and } y\le \varphi(x,d)\,,\\
    g(x,y,d) & \text{ otherwise,}
  \end{cases}
\end{align*}
see Figure~\ref{fig1} for an illustration.
We remark that $f$ is defined on all of $\R^3$, with the first three regions
covering the part outside of $\overline O$. 
The monotonicity of $\varphi$ and (\ref{fixedpoint}) imply that
 $\varphi(y,d)\le p$ for $y\ge p$, therefore the first region is contained in $\{x\le y\}$ and the third in $\{y\le x\}$.
We now show that the fourth region coincides with the set $V$.
Let $(x,y,d)$ be in the fourth region, assume for definiteness $y\ge x$. Since the point is not in the second region, $y> p$. Since it is not in the first one, $x>\varphi(y,d)>|d|/y$. With (\ref{eqfourthregion}) we conclude $(x,y,d)\in V$.

The function $f$ is continuous by definition.
We show that $f\in C^1(\R^3)$. In each of the four regions
which are introduced in the definition of $f$, the 
function  $f$ is smooth, 
and we only
have to consider the partial derivatives along the boundaries. 
In the first region we compute
\begin{align*}
  \partial_y f = \partial_x g \partial_y \varphi + \partial_y g\,,
\end{align*}
but since $\partial_xg(\varphi(y,d),y,d)=0$ by definition of $\varphi$ this
equals $\partial_y g(\varphi(y,d),y,d)$. Hence $\partial_y f$ is
 continuous on the set
$x=\varphi(y,d)$ which defined the boundary between the first and the fourth
region. 
The same holds for $\partial_d f$. The derivative
$\partial_x f$ vanishes in the first region, as well as on the 
boundary between the first and the fourth region. The analogous arguments
show the continuity of $Df$.

The same computation shows that $\partial_yf$, $\partial_x f\ge 0$
everywhere. Therefore $f=h$ on $\Uset$.

To verify the convexity of $f$ we finally 
 compute the second derivatives. Let $\psi=(\varphi(y,d),y,d)$, so that
$f=g\circ \psi$ in the first region. Then
\begin{align*}
  Df = Dg\circ\psi D\psi
\end{align*}
and
\begin{align*}
  D^2f = D^2g \circ\psi D\psi\otimes D\psi + Dg\circ\psi D^2\psi\,.
\end{align*}
From the definition of $\psi$, $D^2\psi_2=D^2\psi_3=0$. From the definition of
$\varphi$, 
$\partial_x g\circ\psi=0$. Therefore the second term vanishes. 

The first term is positive definite because $D^2g\ge 0$ on the set $V$ where
$\partial_xg$, $\partial_yg\ge0$ (see (\ref{beginningchanges})). Therefore $D^2f\ge 0$ in the first
region and by symmetry in the third.
In the second region $D^2f=0$, and in the
fourth $D^2f=D^2g\ge 0$. Therefore $f$ is convex.

\textit{Step~4: $W\qc=W\rc$.} 
From general theory we know that extended-valued polyconvex functions are
quasiconvex, hence $W\pc\le W\qc$. The inequality $W\qc\le W\rc$ holds,
however, only under additional assumptions; we prove it here for the case of
interest. In particular, 
it was shown in \cite{Conti2008} that any quasiconvex function which is finite
on the set of matrices with determinant $t$, for some $t\ne 0$, is rank-one
convex on the same set. Let $t\ne 0$. If $\hdet(t)<\infty$, then $W$ is finite
on $\Sigma_t=\{F: \det F=t\}$. Therefore $W\qc$ is also finite on $\Sigma_t$,
and rank-one convex. The argument of Step  2 only involves laminates within
$\Sigma_t$, therefore we obtain $W\qc\le \Wt$ on $\Sigma_t$. If 
 $\hdet(t)=\infty$ there is nothing to prove.

The case $t=0$ requires a separate treatment. We fix a singular matrix $F$ and
 write $F=a\otimes b$, with $a\ne 0$. The constructions performed in Step 2 
are, in this case, scalar, in the sense that the laminates involve matrices of
the form $F_\pm=a\otimes b_\pm$, for some
other vectors $b_\pm$ and the same $a$ (if $Fv=0$ then one  takes
$F_t=F + t a\otimes w$ instead of $w\otimes w$).
One can then construct a  
scalar test function $\varphi:(0,1)^2\to\R$ with $D\varphi\in\{b_-, b_+\}$ on
large parts of the domain and $\varphi(x)=b\cdot x$ on the boundary, and use
$a\varphi$ in the definition of quasiconvexity (notice that $\det
(Da\varphi)=\det(a\otimes D\varphi)=0$) to conclude.
The proof is now complete.

\section{Discussion}\label{discussion}

\subsection{A formula for the phase boundary}
\label{secformulaphasebdry}
The function $(y,d)\mapsto x=\varphi(y,d)$ 
can be given explicitly in terms of the solutions  of a 
fourth-order polynomial equation in $x^2$. To see this, 
set $L=\lambda^2+\lambda^{-2}$, $M=\lambda^2-\lambda^{-2}$. The 
equation~\bref{eqdxada} reduces to
\begin{align*}
L +M \frac{y^2}{\sqrt{x^2y^2-d^2} }
=\sqrt{2L (x^2+y^2) +
  4 M\sqrt{x^2y^2-d^2} + 8 d}\,.
\end{align*}
Multiplication by $z=\sqrt{x^2 y^2-d^2}$ leads to
\begin{align*}
L z +My^2 =\sqrt{G z^2 + 4 Mz^3 }\,,
\end{align*}
where $G=2L (x^2+y^2) +8d$. We take the square in this identity, find
\begin{align*}
L^2 z^2 +M^2y^4 + 2 LMzy^2 =G z^2+  4 Mz^3 
\end{align*}
and collect terms linear in $z$ on the right-hand side,
\begin{align*}
L^2 z^2 +M^2y^4  -G z^2=  (4 Mz^2 - 2 LMy^2)z\,.
\end{align*}
We square again
\begin{align*}
(L^2 z^2 +M^2y^4  -G z^2)^2=  (4 Mz^2 - 2 LMy^2)^2z^2
\end{align*}
and obtain a rational expression. The terms $G$ and $z^2$ are linear in
$x^2$, the other terms do not depend on $x$. Therefore we have a fourth-order
polynomial equation in $x^2$ for which a solution formula exists.

\begin{figure}
\centering
 \includegraphics[height = .4\textwidth]{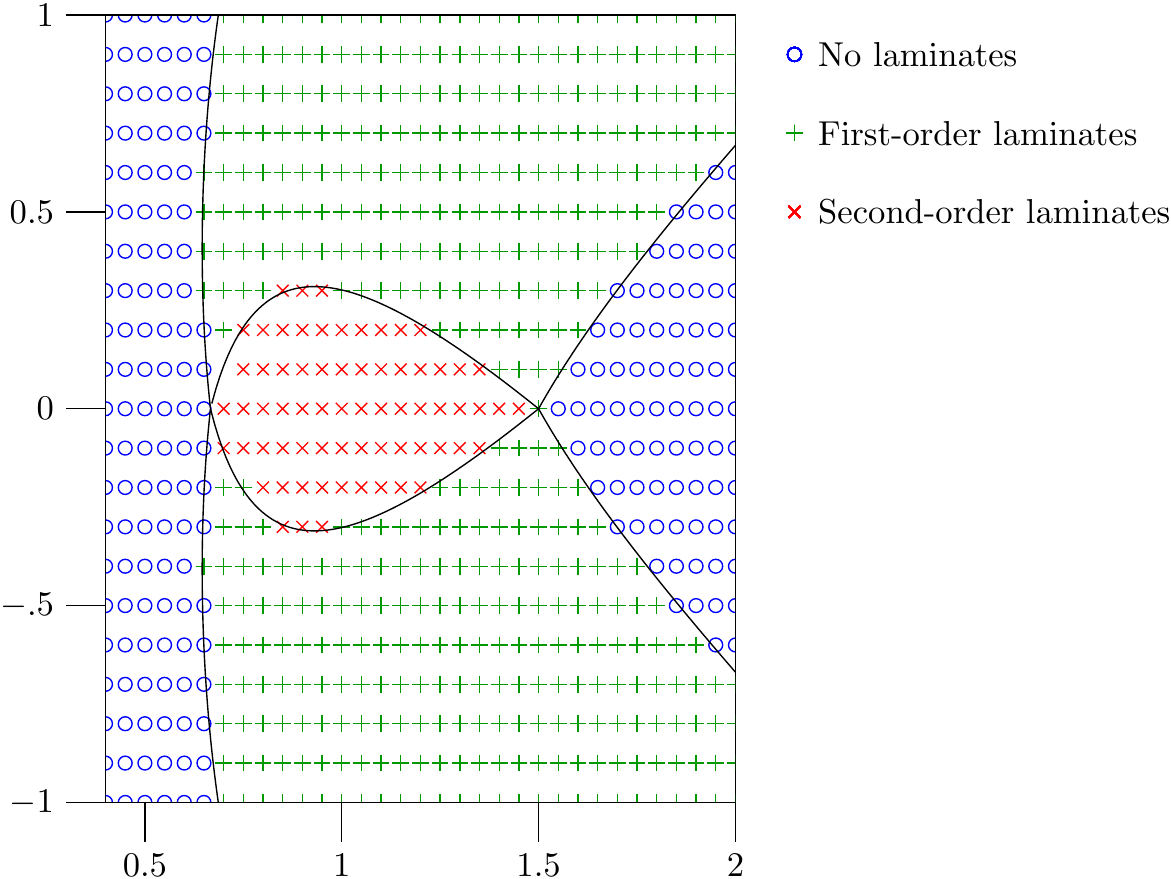}
\caption{\label{fig:phases}
Example of a numerical computation for the relaxed energy using the algorithm
from \cite{ContiDolzmannBaum} for $d=1$, $\lambda=1.5$, with variables as
specified in (\ref{eqFab}). 
The curves are the phase boundaries obtained from the
analytical relaxation, the dots represent numerical results, color-coded to
distinguish the order of lamination. 
The central region is
$K\qc$ and is delimited 
by the  two curves $(x\pm y)^2 + \frac{1}{x^2} = \lambda^2 +
\frac{1}{\lambda^2}$ and corresponds to second-order laminates. The  top and bottom 
regions  correspond to first-order laminates, in the left and right 
regions  the relaxation coincides with the original energy.}
\end{figure}

The plot in Figure~\ref{fig:phases} contains for $d=1$ the zero set of the 
energy density which corresponds to the quasiconvex hull of the two
martensitic wells which is given by~\cite{BallJames92}
\begin{align*}
 K\qc  = \Bset{
F\in \R^{2\times 2}, &\,\det F=1,\,  
|F(e_1\pm e_2)|^2 \leq \lambda^2+\frac{1}{\lambda^2}
}
\end{align*}
or, equivalently, in the coordinates introduced before, 
\begin{align*}
K\qc =\Bset{F\in \R^{2\times 2}, &\,\det F=1,\,x,y
\leq \sqrt{\frac{1}{2}\barg{\lambda^2+\frac{1}{\lambda^2}} }
= \sqrt{\frac{L}{2} }}
\,.
\end{align*}
A short calculation shows that $x=y=\sqrt{L/2}$, $d=1$,
is a solution of~\bref{eqdxada}. This is in agreement with the
results in~\cite{BallJames92} that the relative interior of $K\qc$
is obtained from second order laminates.

\subsection{Comparison with the Ericksen-James energy in two dimensions}
The following energy density originates in~\cite{Ericksen1986,Ericksen1988}
and is usually referred to as the two-dimensional version of the 
 Ericksen-James energy. It is given by
\begin{align*}
 \phi(F) = \kappa_1 \barg{\tr\, C -2 }^2
+ \kappa_2 C_{12}^2 + 
\kappa_3\Barg{ \frac{(C_{11}-C_{22})^2}{4}-\epsilon^2 }^2
\end{align*}
where $C=(C_{ij})_{ij} = F^T F$ and $\epsilon>0$, $\kappa_i>0$, $i=1,2,3$ are
parameters.  
Note that $\phi$ is invariant under the full orthogonal group $\O(2)$ and not only under
$\SO(2)$,
so that the zero set of $\phi$ is given by 
\begin{align*}
 K = \O(2) 
\left( 
\begin{array}{cc} \sqrt{1-\epsilon} & 0 \\ 0 & \sqrt{1+\epsilon}
\end{array}\right)
\cup  \O(2) \left( 
\begin{array}{cc} \sqrt{1+\epsilon} & 0 \\ 0 & \sqrt{1-\epsilon}
\end{array}\right)\,.
\end{align*}
As a consequence, the zero set of the relaxation is much larger than $\SO(2)$.
The relaxation is only known in the special case $\kappa_3=0$ which 
leads to an energy which is convex in the right Cauchy-Green
tensor~\cite{LeDretRaoult1995}. We believe that the relaxation result presented
in Theorem~\ref{maintheorem} will be useful in the design of numerical 
schemes since it provides the full relaxation of a frame-indifferent 
energy which can serve as a model energy for a cubic to tetragonal phase 
transformation in two dimensions.

\subsection{Benchmark example for numerical simulations}
The relaxed energy in Theorem~\ref{maintheorem} can serve as a benchmark
example for numerical schemes for the computation of relaxed energies 
and quasiconvex hulls of sets. In fact, the simulation in Figure~\ref{fig:phases}
was obtained with the algorithm proposed in \cite{ContiDolzmannBaum}. The
relaxation is illustrated for matrices of the form 
\begin{align}\label{eqFab}
 F =
 \begin{pmatrix}
 a & b \\ 0 & 1/a
 \end{pmatrix}
\quad \text{ with }a\in [0.4, 2],\,b\in [-1,1]\,.
\end{align}
The figure shows that the numerically determined phase diagram is in excellent
agreement with the present analytical results. This includes in particular 
 the zero-set of the relaxed energy which is given by
\begin{align*}
 K\qc = \bset{ F\in \R^{2\times 2},\, |F(e_1\pm e_2)|^2 \leq \lambda^2 
+ \frac{1}{\lambda^2}\,,\det F=1 }\,.
\end{align*}
Indeed, the derivation of the analytical result in Theorem~\ref{maintheorem}
was guided by the numerical results in~\cite{ContiDolzmannBaum}.

\subsection{Monotonicity of $\varphi$}
\label{secmonotonicity}
We show here that $\varphi(\cdot,d)$ is decreasing.
From the equation
$\partial_x g(\varphi(y,d),y,d)=0$ and the implicit function theorem we obtain
\begin{equation}
  \partial_x^2 g \partial_y \varphi + \partial_x\partial_y g =0\,.
\end{equation}
Since $\partial_x^2g>0$ where $\partial_xg=0$, it suffices to show that 
$ \partial_x\partial_y g >0$.

From the expressions in Step 3 above we write, using as before $L=\lambda^2+\lambda^{-2}$ and $M=\lambda^2-\lambda^{-2}$, 
\begin{equation*}
  \partial_x\partial_y g 
= \frac{\partial_x A\partial_yA}{2A^{3/2}} + \frac{Mxy}{z^3 A^{1/2}} (2d^2  - x^2y^2)\,.
\end{equation*}
Dropping the positive factor $xy/2A^{3/2}$, it suffices to show that
\begin{equation*}
\Xi=  \frac{\partial_xA}{x}  \frac{\partial_yA}{y} + 2 AM \frac{2d^2 -x^2y^2}{z^3}
\end{equation*}
is nonnegative. Inserting the expressions above and estimating $2d^2 -x^2y^2=d^2-z^2\ge -z^2$
we obtain
\begin{alignat*}1
\Xi&\ge  (L+M\frac{y^2}{z}) (L+M\frac{x^2}{z}) - ( LM(x^2+y^2) + 2 M^2 z + 4 Md) \frac{1}{z}\\
&= L^2 + LM \frac{x^2+y^2}{z} + M^2 \frac{x^2y^2}{z^2} -
LM\frac{x^2+y^2}{z} -2 M^2 - 4 M \frac{d}{z}\\
&= L^2 + M^2 \frac{d^2-z^2}{z^2}- 4 M \frac{d}{z}
\end{alignat*}
where we used $x^2y^2-2z^2=d^2-z^2$.
Since $L^2=M^2+4$ this gives
\begin{alignat*}1
  \Xi&\ge 4 + M^2 \frac{d^2}{z^2}- 4 M \frac{d}{z}\ge0
\end{alignat*}
and concludes the proof.


\end{document}